\documentclass[12pt]{article}
\usepackage{arxiv}
\usepackage[utf8]{inputenc}
\usepackage{amsmath,amssymb,amsfonts}
\title{Extension of Dynamic Mode Decomposition for dynamic systems with incomplete information based on t-model of optimal prediction}
\author{Aleksandr Katrutsa\thanks{Corresponding author} \\
Skolkovo Institute of Science and Technology\\ Moscow, Russia\\
\texttt{aleksandr.katrutsa@phystech.edu} \And Sergey Utyuzhnikov\\
The University of Manchester \\ Manchester, M13 9PL, UK \\
              Moscow Institute of Physics \& Technology\\ Dolgoprudny, Russia\\
 \texttt{s.utyuzhnikov@manchester.ac.uk}
\And 
Ivan Oseledets \\
Skolkovo Institute of Science and Technology\\ Moscow, Russia\\
\texttt{i.oseledets@skoltech.ru} 
}
\date{}

\usepackage{graphicx}
\usepackage{subfig}
\usepackage{algorithm}
\usepackage{algorithmic}

\newcommand{\bx}{\mathbf{x}}
\newcommand{\bu}{\mathbf{u}}

\newcommand{\bz}{\mathbf{z}}
\newcommand{\by}{\mathbf{y}}

\newcommand{\bv}{\mathbf{v}}
\newcommand{\bn}{\mathbf{n}}
\newcommand{\bA}{\mathbf{A}}

\newcommand{\bI}{\mathbf{I}}

\newcommand{\bM}{\mathbf{M}}
\newcommand{\bN}{\mathbf{N}}
\newcommand{\bF}{\mathbf{F}}

\newcommand{\bT}{\mathbf{T}}

\newcommand{\bV}{\mathbf{V}}

\newcommand{\bLambda}{\boldsymbol{\Lambda}}
\newcommand{\blambda}{\boldsymbol{\lambda}}
\newcommand{\bomega}{\boldsymbol{\omega}}

\newcommand{\bX}{\mathbf{X}}
\newcommand{\bbR}{\mathbb{R}}

\newcommand{\calP}{\mathcal{P}}
\newcommand{\calM}{\mathcal{M}}
\newcommand{\calQ}{\mathcal{Q}}

\newcommand{\calI}{\mathcal{I}}
\newcommand{\calL}{\mathcal{L}}

\DeclareMathOperator*{\argmin}{\arg\,\min}

\usepackage{secdot}
\usepackage{xparse}
\usepackage{amsthm} 

\graphicspath{{./fig/}}

\makeatletter
\let\latexparagraph\paragraph
\RenewDocumentCommand{\paragraph}{som}{%
  \IfBooleanTF{#1}
    {\latexparagraph*{#3}}
    {\IfNoValueTF{#2}
       {\latexparagraph{\maybe@addperiod{#3}}}
       {\latexparagraph[#2]{\maybe@addperiod{#3}}}%
  }%
}
\newcommand{\maybe@addperiod}[1]{%
  #1\@addpunct{.}%
}
\makeatother 

\begin{document}

\maketitle

\begin{abstract}
    The Dynamic Mode Decomposition has proved to be a very efficient technique to study dynamic data. 
    This is entirely a data-driven approach that extracts all necessary information from data snapshots which are commonly supposed to be sampled from measurement. 
    The application of this approach becomes problematic if the available data is incomplete because some dimensions of smaller scale either missing or unmeasured.
    Such setting occurs very often in modeling complex dynamical systems such as power grids, in particular with reduced-order modeling. 
    To take into account the effect of unresolved variables the optimal prediction approach based on the Mori-Zwanzig formalism can be applied to obtain the most expected prediction under existing uncertainties. 
    This effectively leads to the development of a time-predictive model accounting for the impact of missing data. 
    In the present paper we provide a detailed derivation of the considered method from the Liouville equation and finalize it with the optimization problem that defines the optimal transition operator corresponding to the observed data.
     In contrast to the existing approach, we consider a first-order approximation of the Mori-Zwanzig decomposition, state the corresponding optimization problem and solve it with the gradient-based optimization method.
     The gradient of the obtained objective function is computed precisely through the automatic differentiation technique.
     The numerical experiments illustrate that the considered approach gives practically the same dynamics as the exact Mori-Zwanzig decomposition, but is less computationally intensive.
\end{abstract}

\section{Introduction}

The data-driven analysis of dynamical systems attracts much attention recently~\cite{shnitzer2020manifold,long2018hybridnet,raissi2018deep,weinan2017proposal}.
However, it is observed that the arbitrary chosen parametric model can be hard to fit from  the given data~\cite{wang2022and}.
To avoid the issues with training, some prior assumptions have to be made to constraint a set of possible models.
One of the such constraints comes from the Mori-Zwanzig decomposition that takes into account hidden dependence between the observed variables and unobserved ones~\cite{fu2020learning}.
The incorporation of such prior assumption in the model leads to recurrent neural network~\cite{wang2020recurrent}.
The combination of Mori-Zwanzig decomposition and Koopman operator formalism is considered in~\cite{lin2021data}, where the interpretation and relations between these approaches are discussed.
Thus, different combinations of data-driven and model-based approaches are studied.

The present study introduces the combination of the  first-order approximation of the optimal prediction approach~\cite{chorin2002optimal} and the Dynamic Mode Decomposition method (DMD)~\cite{schmid2010dynamic,rowley2009spectral}.
The optimal prediction approach is based on  splitting variables into observed and unobserved ones.
This splitting transforms the original dynamical system into coupled ODEs that take into account the effect of unobserved variables to observed ones.
The optimal prediction is originally developed as an algorithm for identification of the observed dynamics~\cite{chorin2003averaging,bernstein2007optimal,chorin2006prediction} and coarse-grained models construction~\cite{li2015incorporation,hudson2020coarse}.
It has a lot of applications in climate modeling~\cite{falkena2019derivation,ghil2020physics,palmer2019stochastic,franzke2020structure}, molecular dynamics simulations~\cite{wang2019implicit,klippenstein2021introducing,volgin2018coarse}, fluid dynamics~\cite{gouasmi2016towards,parish2016reduced}, etc.

In contrast to the optimal prediction, the DMD method is originally focused on the fitting of the transition operator to the measured data.
This method is related to the theory of Koopman operator~\cite{mauroy2020koopman,koopman1931hamiltonian}, which typically motivates the approximation of Koopman operator eigenfunctions by the dynamic modes~\cite{arbabi2017ergodic}.   
The DMD recovers the dynamic modes from the observed data by solving linear least squares problem with the Singular Value Decomposition (SVD)~\cite{kutz2016dynamic}.
Thus, the applicability of the original method is limited~\cite{williams2015data,proctor2016dynamic}.
However, its generalizations to non-linear dependence between sequentially measured samples also exist~\cite{williams2kernel,askham2018variable,degennaro2019scalable}, but they still assume that samples are completely measured.
Therefore, robustness of this approach to noise in the measured data is the topic in many studies~\cite{hemati2017biasing,jovanovic2014sparsity,dawson2016characterizing}.
Thus, the optimal prediction approach and DMD are complementary to each other and the natural idea to combine them is developed in~\cite{curtis2021dynamic}.

The approach suggested in~\cite{curtis2021dynamic} considers the discretizations of the both observed variables equation and memory kernel equation.
After combining the discretized equations, the authors state an optimization problem with respect to the transition operator that is naturally tuned to the related operator from the DMD. 
To solve the stated optimization problem, the approximation of the first-order optimality condition is derived and the approximate solution is written in a closed form.
Thus, although the original model is exact, the solution is computed only approximately without detailed description of the intermediate steps. 
We have noted that instead of considering the exact model, its approximation with the t-model~\cite{chertock2008modified} can be used.
This approach significantly simplifies further derivation of the optimization problem and the optimization problem itself.
Then, we solve it accurately with the adaptive  gradient-based method called  Adam~\cite{kingma2014adam}.
The gradient of the objective function is computed with an automatic differentiation technique~\cite{baydin2018automatic}.
In the numerical experiment section, we demonstrate that the proposed approach leads to practically the same dynamics of the test model, more robust to the noise and requires essentially less intensive computations to get the optimal transition operator.
To compare the proposed approach with the optimal prediction, we correct the original equations from~\cite{curtis2021dynamic} and provide the detailed derivation of the final optimization problem for reader convenience.
The developed method allows us to control the approximation accuracy and the computational efficiency simultaneously.

The rest of the paper is organized as follows. 
In Section~\ref{sec::mz_decomposition}, the Mori-Zwanzig formalism is applied for the analysis of dynamic data assigned with resolvable (measured) variables. 
Then, a variational problem for the transition operator is obtained that effectively leads to a finite dimensional approximation of the Koopman operator. 
Thanks to the Mori-Zwanzig decomposition, it contains a memory term that is responsible for the effect of unresolved variables. 
Without the memory term the transition operator coincides with the standard DMD operator. 
In turn, a more accurate transition operator must be the memory dependent DMD as first noted in~\cite{curtis2021dynamic}. 
This kind of approximation is simplified in Section~\ref{sec::t-model} with the use of the modified optimal prediction method that represents a linear expansion of the Mori-Zwanzig decomposition operator with respect to time. 
The efficiency of the developed algorithm is tested in Section~\ref{sec::num_exp} on the model of a coupled oscillator. 
This section is followed by the conclusion in Section~\ref{sec::conclusion}.

\section{Mori-Zwanzig decomposition as a nonlinear extension of Dynamic Mode Decomposition}
\label{sec::mz_decomposition}
The discussed above DMD approach relies on the linear dependence between sequential data samples: $\bx_{i+1} \approx \bA_{dmd}\bx_i$, where $i=1,\ldots,m-1$.
This approach reduces to the solving of the following linear least-squares problem with respect to the transition matrix~$\bA_{dmd}$
\begin{equation}
    \bA_{dmd} = \argmin_{\bA}\| \bX_+ - \bA \bX_-\|^2_F,
    \label{eq::dmd_opt_problem}
\end{equation}
where $\bX_+ = [\bx_m, \ldots, \bx_2]$, $\bX_- = [\bx_{m-1}, \ldots, \bx_1]$ are stacked samples $\bx_i \in \mathbb{R}^n$ in the matrices of size $n \times (m-1)$.
Problem~\eqref{eq::dmd_opt_problem} has a well-known analytical solution in the form $\bA_{dmd} = \bX_+\bX^{\dagger}_-$, where $\bX^{\dagger}_-$ is the Moore-Penrose pseudoinverse matrix~\cite{kutz2016dynamic}.
However, complex real-world data is usually governed by the non-linear dynamical systems and cannot be approximated by linear models with sufficiently high accuracy.
To address this issue, a hybrid approach can be used.
This algorithm is based on the assumption that we observe only a part of variables describing the considered system.
Therefore, we should take into account not only the dynamics of the observed variables, but also its interaction with the unobserved variables.
This assumption leads to a more complicated problem than~\eqref{eq::dmd_opt_problem} and in further section this problem is derived and discussed in detail.



\subsection{Mori-Zwanzig decomposition}
Consider the Liouville equation with operator $\mathcal{L} = \sum_{i=1}^n f_j(\bz) \frac{\partial }{\partial z_j}$ in the right-hand side and initial value $\bu_0$:
\begin{equation}
    \begin{aligned}
    & \frac{\partial\bu}{\partial t} = \mathcal{L}\bu,\\
    & \bu(\bz, 0) = \bu_0,
    \label{eq::liouv_general}
\end{aligned}
\end{equation}
where $\bu \in \bbR^n$.
Eq.~\eqref{eq::liouv_general} is equivalent to the Cauchy problem for its characteristic equation:
\begin{equation}
\begin{split}
    &\frac{d\mathbf{\Phi}}{dt} = \mathbf{f}(\mathbf{\Phi}), \quad \mathbf{\Phi} \in \bbR^n,\\
    &\mathbf{\Phi}(0) = \bz.
\end{split}
\label{eq::caushy_char}
\end{equation}
Indeed, it immediately follows that if $\bu_0 = \bz$, then $\bu(\bz, t) = \mathbf{\Phi}(\bz, t)$. 
A principal difference between these two equations is that Eq.~\eqref{eq::caushy_char} can be nonlinear while Eq.~\eqref{eq::liouv_general} is always linear.

Below, we use so-called semigroup notation and write the solution as $\bu(t) = e^{\calL t}\bu_0$.
Now we split the unknown variables $\bu(t)$ into two parts $\bu(t) = [\bx(t), \by(t)]$, where $\bx(t)$ are the observed variables, also known as \emph{resolved variables}, and $\by(t)$ are the \emph{unresolved variables} that do not available after measurements.
Therefore, we can project~\eqref{eq::liouv_general} onto the space of only resolved variables with projector~$\mathcal{P}$, i.e. $\mathcal{P}\bu(t) = \bx(t)$.
In particular, we can define this projector as $\calP \bu(t) = \mathbb{E} [\bu(t) | \bx(t)]$, where the expectation is computed with respect to the probability distribution of the unresolved variables $\by(t)$.
In addition, denote by $\mathcal{Q}$ the projector onto the space of unresolved variables, i.e. $\mathcal{Q}\bu(t) = \by(t)$.
Note that $\calP + \calQ = \calI$, which is identity operator.
After that we can re-write~\eqref{eq::liouv_general} only for the resolved variables $\bx(t)$ but take into account the effect from the unresolved variables $\by(t)$:

\begin{equation}
\begin{split}
    \frac{\partial(\mathcal{P}e^{\mathcal{L}t}\bu_0)}{\partial t} &= \calP \calL e^{\calL t}\bu_0 = \calP e^{\calL t} \calI \calL \bu_0 = \calP e^{\calL t} (\calP + \calQ) \calL \bu_0 \\ 
    &= \calP e^{\calL t} \calP \calL \bu_0 + \calP e^{\calL t}\calQ\calL \bu_0,
\end{split}
\end{equation}
where we use equality $\calL e^{\calL t} = e^{\calL t} \calL$. 
The first item $\calP e^{\calL t} \calP \calL \bu_0$ indicates the dynamics of the resolved variables.
The second item $\calP e^{\calL t}\calQ\calL \bu_0$ is processed as follows.
From the Dyson operator identity
\begin{equation}
    e^{\calL t} = e^{\calQ \calL t} + \int_{0}^t e^{\calL (t - \tau)}\calP \calL e^{\calQ \calL \tau} d\tau
    \label{eq::dyson}
\end{equation}
it follows that the second term can be re-written in the form:
\begin{equation*}
    \begin{split}
        \calP e^{\calL t}\calQ\calL \bu_0 &= \calP e^{\calQ \calL t}\calQ \calL \bu_0 + \calP\int_{0}^t e^{\calL (t - \tau)}\calP \calL e^{\calQ \calL \tau} \calQ \calL \bu_0 d\tau \\
        &= \calP \by(t) + \int_{0}^t \calP e^{\calL (t - \tau)} \calM (\tau, \bu_0) d\tau,
    \end{split}
\end{equation*}
where $\calM (\tau, \bu_0) = \calP \calL e^{\calQ \calL \tau} \calQ \calL \bu_0 = \calP \calL \bN(\tau, \bu_0)$.
By $\bN(\tau, \bu_0)$ we denote the dynamics $e^{\calQ \calL \tau} \calQ \calL \bu_0$ in the space of unresolved variables that can be interpreted as noise~\cite{curtis2021dynamic}.
Since subspace of resolved variables $\bx(t)$ is orthogonal to the subspace of unresolved ones, $\calP \by(t) \equiv 0$.
Thus, we get the initial value problem for the resolved variables with a memory term depending on the unresolved variables that proceeds further:
\begin{equation}
    \frac{\partial \bx}{\partial t} = \calP \calL \bx + \int_0^{t} \calP e^{\calL (t - \tau)} \calM(\tau, \bu_0) d\tau,
    \label{eq::mz_1}
\end{equation}
where $\calM (\tau, \bu_0) = \calP \calL e^{\calQ \calL \tau} \calQ \calL \bu_0$ and $\calP e^{\calL t} \calP \calL \bu_0 = \calP \calL \calP e^{\calL t} \bu_0 =\calP \calL \bx$.
Thus, the second term in~\eqref{eq::mz_1} represents the effect of unresolved variables on resolved ones.
It is clear that this effect is nonlocal and not determenistic.

To exclude the explicit action of operator $e^{\calQ \calL \tau}$ in the space of unresolved variables, we need to derive the equation for kernel $\calM (\tau, \bu_0)$.
From the Dyson operator identity~\eqref{eq::dyson} we can get the equation for $\bN(\tau, \bu_0)$:
\begin{equation*}
    \begin{split}
        \bN(\tau, \bu_0) &= e^{\calQ \calL \tau} \calQ \calL \bu_0 = e^{\calL \tau} \calQ \calL \bu_0 - \int_0^{\tau} e^{\calL (\tau - s)}\calP \calL e^{\calQ \calL s} \calQ \calL \bu_0 ds \\
        &= e^{\calL \tau} \bN(0, \bu_0) - \int_0^{\tau} e^{\calL (\tau - s)} \calM (s, \bu_0)ds
    \end{split}
\end{equation*}
and after multiplication by $\calP\calL$ both sides of this equation and thanks to commutability of $\calL$ and $e^{\calL \tau}$, we have the equation for kernel $\calM(\tau, \bu_0)$:
\begin{equation}
    \calM(\tau, \bu_0) + \int_0^{\tau} \calP e^{\calL (\tau - s)} \calL \calM (s, \bu_0)ds = \calP e^{\calL \tau} \calL \bN(0, \bu_0).
    \label{eq::mz_2}
\end{equation}

Thus, we have equations~\eqref{eq::mz_1} and~\eqref{eq::mz_2} that completely define the dynamics of resolved variables and its interaction with the unresolved variables.
In the next paragraph we consider the discretization of these equations.
From the discretized equations, an optimization problem follows to fit the unknown transition operator for resolved variables.
The main feature of this operator is that it  incorporates interactions of resolved variables with unresolved ones.

\paragraph{Discretization of the Mori-Zwanzig decomposition. } 

To derive the proper discretization, we firstly introduce notation for the target operator.
Denote by $\bT_{mz}(t) = \calP e^{\calL t}$ the projection of the evolution operator in the resolved variables subspace.
Then, the diagonalized form of this operator in finite dimension looks like
\[
\bT_{mz}(t) \approx \bV e^{\bLambda t}\bV^{-1},
\]
which means that $\calP\calL \approx \bV\bLambda \bV^{-1} = \bA_{mz}$, where $\bLambda$ is a diagonal matrix.

Now, following the procedure in~\cite{curtis2021dynamic} we use the first-order approximation of the left-hand side and the trapezoidal method to approximate the integral in the right-hand side of~\eqref{eq::mz_1}.
After that the discretization of~\eqref{eq::mz_1} has the following form
\begin{equation}
\begin{split}
     \frac{\bx_{n+1} - \bx_n}{\Delta \tau} &= \calP\calL \bx_n + \Delta \tau \left( \sum_{k=1}^{n-1} \calP e^{(n - k)\Delta\tau\calL} \calM(k\Delta\tau, \bu_0) + \right.\\ &\left.\frac{1}{2}[\calP \calM(n\Delta\tau, \bu_0) + \calP e^{\calL n\tau}\calM(0, \bu_0)] \right).
\end{split}
    \label{eq::mz_discr_1}
\end{equation}
Here, vectors $\bx_n$ coincide with the observed samples that are columns of matrices $\bX_+$ and $\bX_-$ from~\eqref{eq::dmd_opt_problem}.
Now we reduce this discretized equation to the form such that we can state a similar optimization problem for the transition operator $\bA_{mz}$.

Next, denote by $\calM_k$ the vector $\calM(k\Delta \tau, \bu_0)$ and re-write~\eqref{eq::mz_discr_1} in the basis of eigenvectors $\bV$.
In such transformation we use the following notations: $\hat{\bx}_n = \bV^{-1}\bx_n$ and $\widehat{\calM}_k = \bV^{-1}\calM_k$.
Then we arrive at the following equation for $\hat{\bx}_n$:
\begin{equation}
    \hat{\bx}_{n+1} = (\bI + \Delta \tau\bLambda)\hat{\bx}_n + \frac{\Delta \tau^2}{2} \left(\widehat{\calM}_n + e^{n\Delta\tau\bLambda}\widehat{\calM}_0 + 2 \sum_{k=1}^{n-1}e^{(n-k)\Delta \tau \bLambda}\widehat{\calM}_k \right).
    \label{eq::x_discr}
\end{equation}
Now, we need to discretize~\eqref{eq::mz_2}, re-write it in the same basis $\bV$ and get the equation for $\widehat{\calM}_k$. 
But before discretization, let us re-write~\eqref{eq::mz_2} in the following form:
\begin{equation}
\begin{split}
    &\calM(\tau, \bu_0) + \int_0^{\tau} \calP e^{\calL (\tau - s)} (\calP + \calQ)\calL \calM (s, \bu_0)ds \\
    &= \calM(\tau, \bu_0) + \int_0^{\tau} \calP e^{\calL (\tau - s)}\calP \calL \calM (s, \bu_0)ds = \calP e^{\calL \tau} \calL \bN(0, \bu_0)
\end{split}
    \label{eq::reduced_mz_2}
\end{equation}
since $\calI = \calP + \calQ$ and $\int_0^{\tau} \calP e^{\calL (\tau - s)}\calQ\calL \calM (s, \bu_0)ds = 0$.
To discretize Eq.~\eqref{eq::reduced_mz_2}, we again use the trapezoidal method to approximate the integral in the left-hand side and use approximation $\calL \bN(0, \bu_0) \approx \calM(0, \bu_0)$:
\begin{equation*}
    \begin{split}
        &\calM_n + \frac{\Delta \tau}{2} \left( \calP\calL \calM_n + \calP e^{n\Delta\tau\calL}\calP\calL \calM_0 + 2 \sum_{k=1}^{n-1} \calP e^{\calL (n - k)\Delta\tau}\calP\calL \calM_k \right) = \\
        &= \calP e^{n\Delta\tau \calL}\calM_0
    \end{split}
\end{equation*}
or in basis $\bV$:
\[
\widehat{\calM}_n + \Delta \tau \left( \frac12 (\bLambda \widehat{\calM}_n + e^{\bLambda n\Delta\tau} \bLambda \widehat{\calM}_0) + \sum_{k=1}^{n-1} e^{\bLambda (n-k)\Delta\tau}\bLambda \widehat{\calM}_k \right) = e^{\bLambda n\Delta\tau} \widehat{\calM}_0.
\]
From the above equation, it follows that 
\begin{equation}
    \widehat{\calM}_n = \left( \bI + \frac{\Delta\tau}{2}\bLambda \right)^{-1} \left( e^{\bLambda n\Delta\tau} \left(\bI - \frac{\Delta\tau}{2}\bLambda \right)\widehat{\calM}_0 - \Delta\tau \bLambda  \sum_{k=1}^{n-1} e^{\bLambda (n-k)\Delta\tau} \widehat{\calM}_k \right).
    \label{eq::mz_memory_discr}
\end{equation}
To insert the equation for memory term $\widehat{\calM}_k$ into the equation for resolved variables~\eqref{eq::x_discr}, one writes the following recurrent equation:
\[
\widehat{\calM}_n = e^{\Delta\tau\bLambda}\bM(\bLambda)\widehat{\calM}_{n-1},\; n \geq 1,
\]
where $\bM(\bLambda) = \bI - \Delta \tau \bLambda \left( \bI + \frac{\Delta\tau}{2}\bLambda\right)^{-1}$.
Thus, equation~\eqref{eq::x_discr} can be re-written as follows:
\begin{equation}
        \hat{\bx}_{n+1} = (\bI + \Delta \tau \bLambda)\hat{\bx}_n + \frac{\Delta \tau^2}{2}e^{n\Delta\tau \bLambda}\left(\bM^n(\bLambda) + \bI + 2\sum_{k=1}^{n-1} \bM^k(\bLambda) \right)\widehat{\calM}_0.
        \label{eq::x_discr_rec}
\end{equation}
Now, consider the expression with matrix $\bM(\bLambda)$ in more detail and simplify it:
\begin{equation*}
    \begin{split}
        & \bM^n(\bLambda) + \bI + 2\sum_{k=1}^{n-1} \bM^k(\bLambda) = \bM^n(\bLambda) - \bI + 2\sum_{k=0}^{n-1} \bM^k(\bLambda)\\ &=\bM^n(\bLambda) - \bI + 2(\bI - \bM(\bLambda))^{-1}(\bI - \bM^n(\bLambda)) \\
        &= 
        (\bM^n(\bLambda) - \bI)(\bI + 2(\bM(\bLambda) - \bI)^{-1})
        =(\bM^n(\bLambda) - \bI)\left(-\frac{2}{\Delta\tau} \bLambda^{-1} \right).
    \end{split}
\end{equation*}
Thus, we simplify equation~\eqref{eq::x_discr_rec} to the following form:
\begin{equation}
    \hat{\bx}_{n+1} = (\bI + \Delta \tau \bLambda)\hat{\bx}_n -  \Delta \tau \bLambda^{-1} e^{n\Delta\tau \bLambda}(\bM^n(\bLambda) - \bI) \widehat{\calM}_0.
    \label{eq::x_discr_lamb}
\end{equation}
Denote by $\bar{\bLambda}$ the diagonal matrix $\bI + \Delta \tau \bLambda$ and by $\bA$ the matrix $\bV\bar{\bLambda}\bV^{-1}$.
Then, since $\bLambda = \frac{\bar{\bLambda} - \bI}{\Delta\tau}$ we can re-write~\eqref{eq::x_discr_lamb} in terms of $\bar{\bLambda}$:
\begin{equation}
    \hat{\bx}_{n+1} = \bar{\bLambda}\hat{\bx}_n -  \Delta \tau^2 (\bar{\bLambda} - \bI)^{-1}e^{n(\bar{\bLambda} - \bI)}(\bM^n(\bar{\bLambda}) - \bI) \widehat{\calM}_0,
    \label{eq::x_discr_lamb_hat}
\end{equation}
where $\bM(\bar{\bLambda}) = \bI - 2 (\bar{\bLambda} - \bI)(\bar{\bLambda} + \bI)^{-1}$ or in terms of the transition operator~$\bA$:
\begin{equation}
    \bx_{n+1} = \bA\bx_n -  \Delta \tau^2(\bA - \bI)^{-1} e^{n(\bA - \bI)}(\bM^n(\bA) - \bI) \calM_0,
    \label{eq::x_discr_A}
\end{equation}
where $\bM(\bA) = \bI - 2 (\bA - \bI)(\bA + \bI)^{-1}$.

Since we get equation~\eqref{eq::x_discr_A} that relates two sequential snapshots $\bx_n$ and $\bx_{n+1}$, we can state the optimization problem to find the transition operator $\bA_{mz}$ corresponding to the Mori-Zwanzig decomposition approach:
\begin{equation}
\bA_{mz} = \argmin_{\bA} \| \bX_+ - \bA \bX_- + \Delta \tau^2 \widetilde{\bM}(\bA, \bn) \|_F^2,
    \label{eq::mzdmd_final_opt}
\end{equation}
where $\widetilde{\bM}(\bA, \bn) = (\bA - \bI)^{-1}\bF(\bA, \bn)$ and 
\[
\bF(\bA, \bn) = 
\begin{bmatrix} \mathbf{0} & F_1(\bA, \bn) & \ldots & F_{m-2}(\bA, \bn) \end{bmatrix}
\]
is a matrix of size $n\times (m-1)$ such that $F_j(\bA, \bn) = e^{j(\bA - \bI)}(\bM^j(\bA) - \bI)\bn$. 
We introduce a vector~$\bn$ that corresponds to the particular initialization of the memory term $\calM_0$.

\section{t-model and DMD}

\label{sec::t-model}

In this section, we consider the t-model~\cite{chertock2008modified} as an approximation of the exact Mori-Zwanzig decomposition~\eqref{eq::mz_1},~\eqref{eq::mz_2}.
We also present the discretization scheme of this model similar to~\eqref{eq::x_discr},~\eqref{eq::mz_memory_discr} corresponding to the exact Mori-Zwanzig decomposition.

To obtain equations for the t-model, we re-write the integral in the right-hand side of~\eqref{eq::mz_1} in the following way
\[
\frac{\partial \bx}{\partial t} = \calP \calL \bx + t \calP e^{\calL t} \calM(0, \bu_0) + \mathcal{O}(t^2)
\]
and truncate the terms of higher order than linear in the right-hand side.
Thus, we get the following approximate equation for the resolved variable dynamics
\begin{equation}
    \frac{\partial \bx}{\partial t} = \calP \calL \bx + t \calP e^{\calL t} \calM(0, \bu_0).
    \label{eq::tmodel}
\end{equation}
In the discretization of this equation we do not need equation~\eqref{eq::mz_memory_discr} for $\calM_k$ that significantly simplifies the optimization problem which is derived further.

The discretization of~\eqref{eq::tmodel} and transformation in the basis $\bV$ give the following equation 
\[
\hat{\bx}_{n+1} = (\bI + \Delta \tau \bLambda)\hat{\bx}_n + n\Delta\tau^2 e^{n\Delta\tau \bLambda} \widehat{\calM}_0,
\]
where we presume equality $t = n\Delta \tau$.
After that, we can re-write this equation in terms of matrix $\bar{\bLambda}_t = \bI + \Delta \tau \bLambda$ similar to~\eqref{eq::x_discr_lamb_hat}:
\[
\hat{\bx}_{n+1} = \bar{\bLambda}_t\hat{\bx}_n + n\Delta\tau^2 e^{n(\bar{\bLambda}_t - \bI)} \widehat{\calM}_0,
\]
and state the final optimization problem for the transition operator $\bA_{t} = \bV\bar{\bLambda}_t\bV^{-1}$ corresponding to the  t-model:
\begin{equation}
    \bA_{t} = \argmin_{\bA} \| \bX_+ - \bA \bX_- - \Delta \tau \widetilde{\bM}(\bA, \bn) \|_F^2,
    \label{eq::final_tmodel_opt}
\end{equation}
where $\widetilde{\bM}(\bA, \bn) = \begin{bmatrix} 0 & g_1(\bA, \bn) & \ldots & g_{m-2}(\bA, \bn)  \end{bmatrix}$ and $g_j(\bA, \bn) = j\Delta \tau e^{j(\bA - \bI)}\bn$.
Comparison of expressions for $\widehat{\bM}$ and $\widetilde{\bM}$ highlights the simplicity of the t-model. 
Further numerical experiments demonstrate that this simplified model is still sufficiently accurate for a test dynamical system used in previous studies~\cite{chorin2002optimal,curtis2021dynamic}.

\section{Numerical experiments}
\label{sec::num_exp}

In this section we present the comparison of the models described above.
We illustrate the main features of the considered methods such as approximation accuracy, running time and robustness to noise.
We also demonstrate the limitations of the presented model to simulate a simple dynamical system. 

\paragraph{Reconstruction of the dynamics from the averaged eigendecomposition}
Since optimization problems~\eqref{eq::mzdmd_final_opt} and~\eqref{eq::final_tmodel_opt} depend on the memory term initialization $\bn$, we have to average the result with respect to the different values of vector $\bn$.
To perform such averaging procedure we sample $N_u$ vectors $\bn$, solve problem~\eqref{eq::mzdmd_final_opt} for every generated $\bn$ and then average the spectra and corresponding eigenvectors over the solutions of these problems. 
Having the averaged spectrum $\bar{\blambda}$ and corresponding eigenvectors in columns of matrix $\bar{\bV}$, we can reconstruct the expectation of of the considered system dynamics.
If the initial state $\bx_0$ is given, the expectation $\bar{\bx}(t)$ at the time $t$ can be computed according to the following equation:
\begin{equation}
    \bar{\bx}(t) = \bar{\bV} \mathrm{diag}(e^{\bomega t}) \bar{\bV}^{-1}\bx_0, \quad \bomega = \frac{\log (\bar{\blambda})}{\Delta \tau},
    \label{eq::dynamic_recontr}
\end{equation}
where the exponent and logarithm functions are elementwise and $\mathrm{diag}(\bx)$ denotes the diagonal matrix with vector $\bx$ in the diagonal.
Since we perform averaging of the derived spectra and corresponding eigenvectors, we need to analyze the variance of this estimate.
To compute the variance $\bv(t)$ we use the empirical estimation formula:
\begin{equation}
    \bv(t) = \frac{1}{N_u}\sum_{i=1}^{N_u} (\bx_i(t) - \bar{\bx}(t))^2,
    \label{eq::variance_estimate}
\end{equation}
where $\bx_i(t)$ is dynamics obtained from solving problem~\eqref{eq::mzdmd_final_opt} with an  $i$-th initialization of vector~$\bn$.
This dynamics can be recovered from the $i$-th solution of the corresponding optimization problem through~\eqref{eq::dynamic_recontr}, where the spectrum and eigenvectors of the solution are used instead of $\bar{\blambda}$ and $\bar{\bV}$.  

\paragraph{Hamiltonian system}
Next, consider the simple model of coupled oscillators from~\cite{chorin2002optimal,curtis2021dynamic}:
\begin{equation}
    \begin{cases}
    \dot{y}_1 = y_2 \\
    \dot{y}_3 = y_4\\
    \dot{y}_2 = -y_1(1 + y_3^2)\\
    \dot{y}_4 = -y_3(1 + y_1^2),
    \end{cases}
    \label{eq::simple_model}
\end{equation}
where $y_1(0) = \hat{x}_1$ and $y_2(0) = \hat{x}_2$ are supposed to be the resolved variables and fixed.
At the same time, $y_3(0)$ and $y_4(0)$ are the unresolved and therefore the  corresponding memory term is modeled with the  Gaussian noise, i.e. $y_3(0) = \tilde{x}_3$ and $y_4(0) = \tilde{x}_4$ such that $\tilde{x}_j \sim \mathcal{N}(0, \sigma^2)$ for $j=3,4$ and $\mathcal{N}(0, \sigma^2)$ denotes Gaussian distribution with zero mean and the  variance $\sigma^2$.  

\begin{figure}[!ht]
    \centering
    \subfloat{\includegraphics[width=0.45\textwidth]{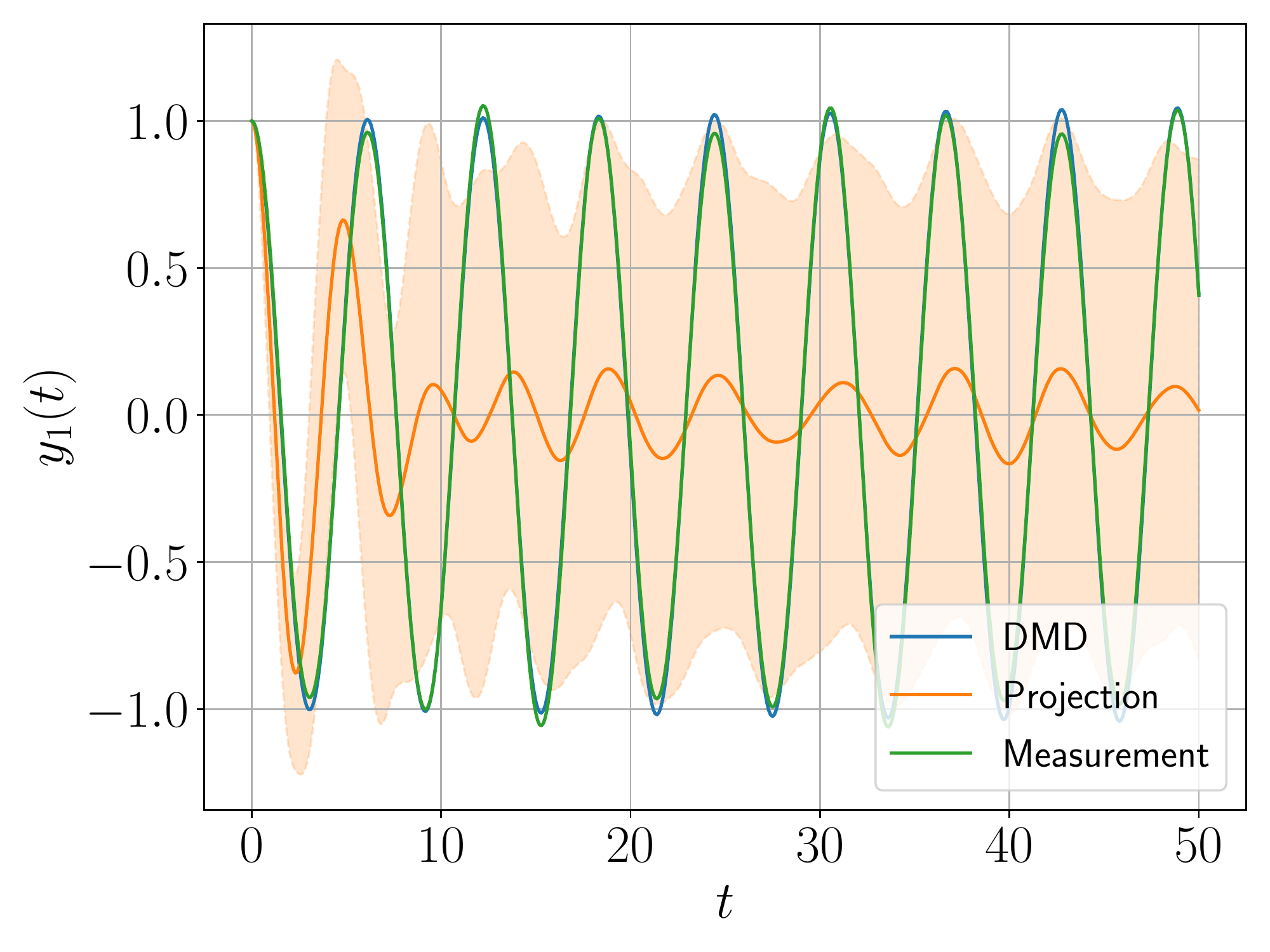}}
    ~
    \subfloat{\includegraphics[width=0.45\textwidth]{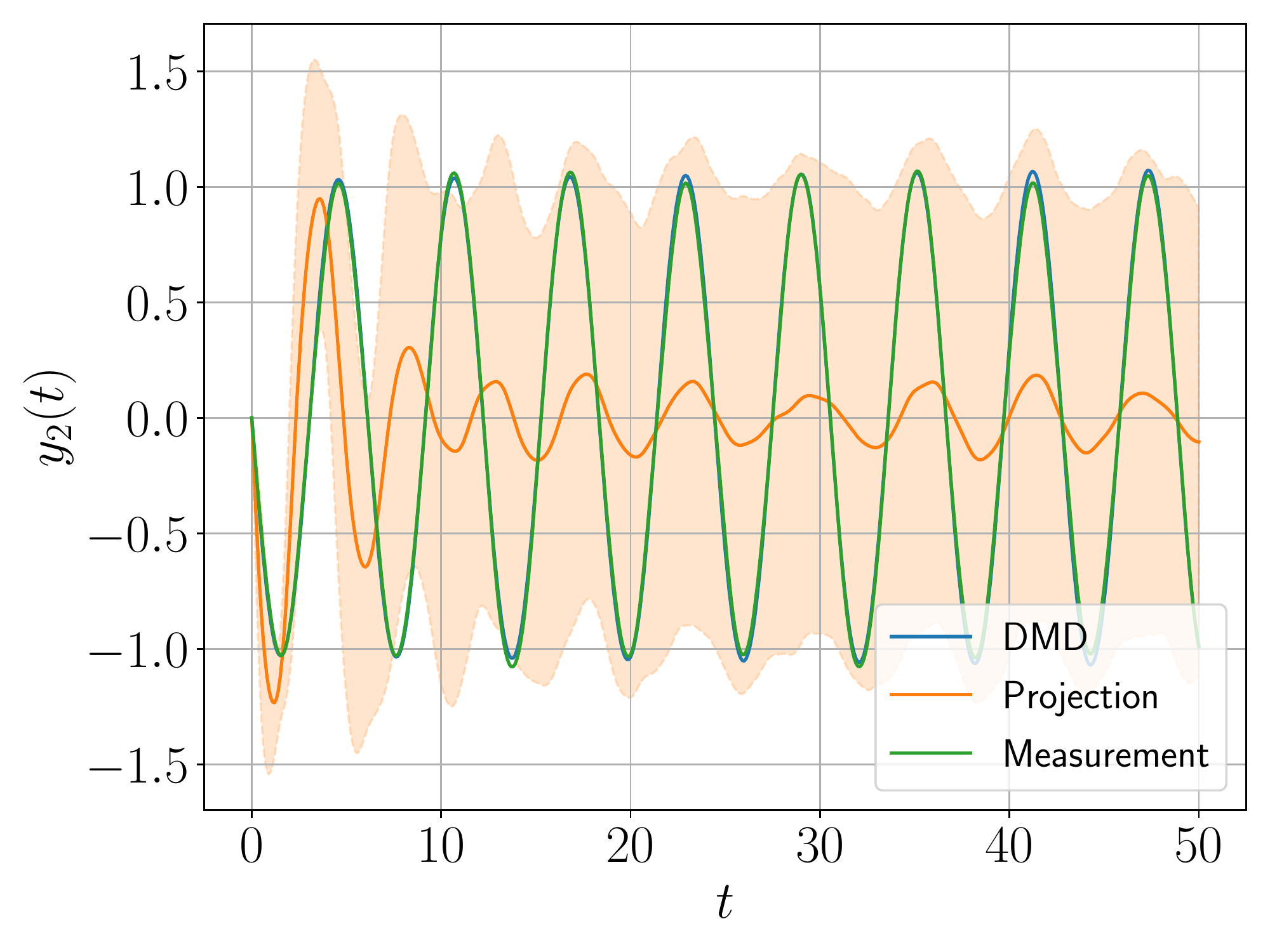}}
    \caption{Comparison of the measured and projected dynamics with DMD reconstruction for resolved $y_1(t)$ and $y_2(t)$. The difference between projected dynamics and DMD illustrates the effect of coupling between resolved and unresolved variables. 
    The filled space indicates the magnitude of the variance corresponding to the Projection approach.}
    \label{fig::dmd_vs_true}
\end{figure}

A naive approach to simulate the dynamics of the system~\eqref{eq::simple_model} is to sample a pair of initial values $(\tilde{x}_3, \tilde{x}_4)$ and run integration~\cite{hindmarsh1983odepack} with these initial values.
We denote the result of such simulation as  \emph{Measurement} and plot it in Figure~\ref{fig::dmd_vs_true}.
Based on this simulation we construct matrices $\bX_+$ and $\bX_-$ from the DMD problem statement and solve  problem~\eqref{eq::dmd_opt_problem}.
Then, we use the obtained operator $\bA_{dmd}$ to build the reconstructed dynamics according to~\eqref{eq::dynamic_recontr}, where the spectrum and eigenvectors of $\bA_{dmd}$ are used instead of averaged ones.
The reconstructed dynamics is plotted in Figure~\ref{fig::dmd_vs_true} and labeled by \emph{DMD}.

To simulate the expectations of the resolved variables, we sample $N=10^3$ values for $(\tilde{x}_3, \tilde{x}_4)$, run the numerical integrator~\cite{hindmarsh1983odepack} and average the results.
The result of this procedure is labeled as \emph{Projection} in Figure~\ref{fig::dmd_vs_true}, where variance estimation of the derived projected dynamics is also shown.
Thus, in contrast to the \emph{DMD} and \emph{Measurement} results, the averaging procedure is performed for many initial values of unresolved variables $(\tilde{x}_3, \tilde{x}_4)$.   
Other parameters are $\Delta \tau = 10^{-1}$, time range $t \in [0, 50]$, number of points in time grid is $N_t = 501$, and $\sigma =1$.
The initial values for the resolved variables  are $(\hat{x}_1, \hat{x}_2) = (1, 0)$.
The obtained dynamics demonstrate a  significant difference between the considered approaches, see Figure~\ref{fig::dmd_vs_true}.


\begin{figure}[!ht]
    \centering
    \subfloat{\includegraphics[width=0.45\textwidth]{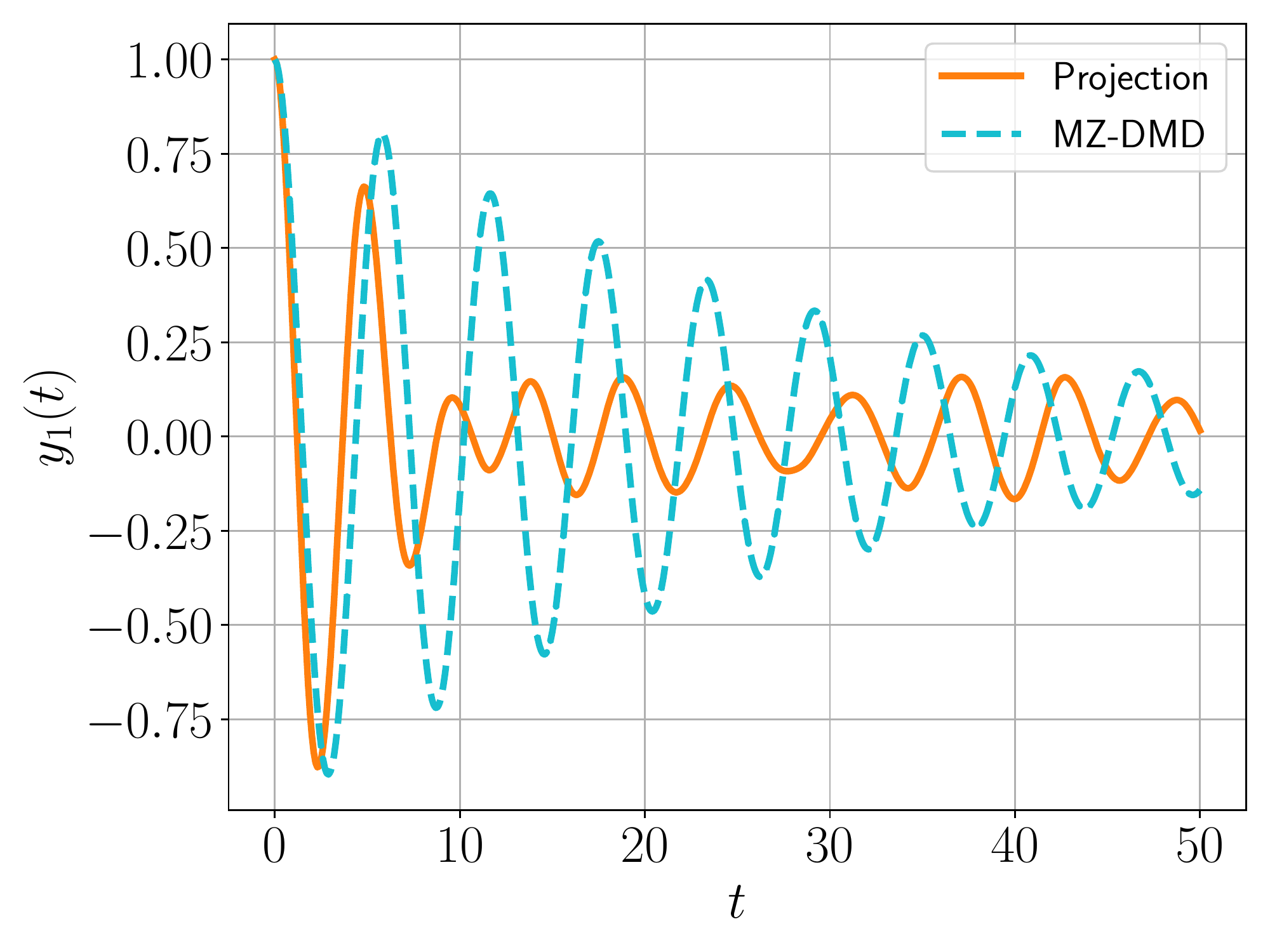}}
    ~
    \subfloat{\includegraphics[width=0.45\textwidth]{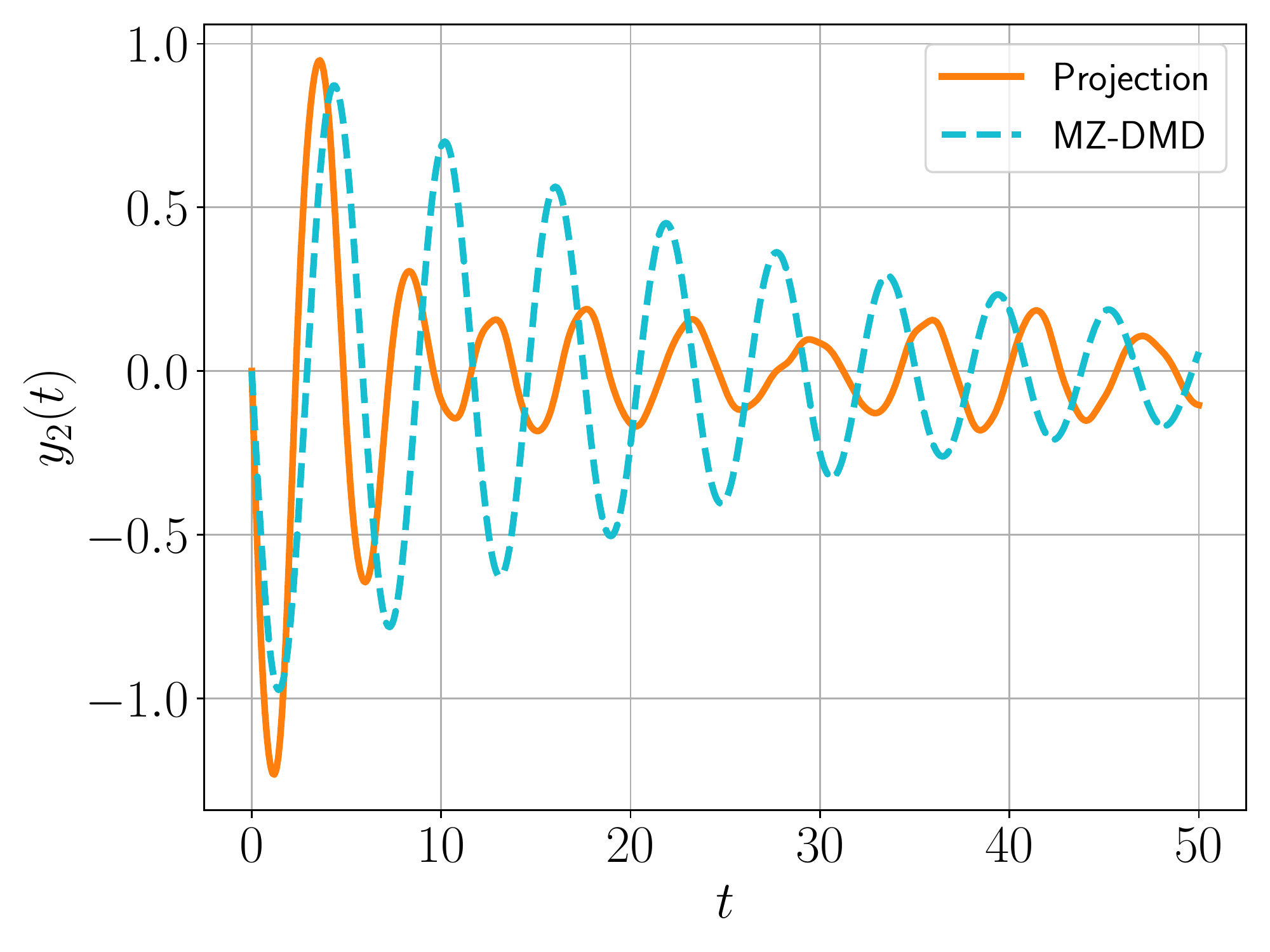}}
    \caption{Comparison of the projected dynamics and MZ-DMD reconstruction for resolved $y_1(t)$ and $y_2(t)$. Although the dynamics do not coincide, the decay trend of the projected dynamics is captured by the dynamics reconstructed from the MZ-DMD approach. We use $N_u = 100$ different initialization of vector~$\bn$. The variance estimate is smaller than $10^{-1}$ uniformly on the considered time range and is not shown in plots. This observation demonstrates the robustness of the reconstructed dynamics with respect to different initialization~$\bn$ of memory term.}
    \label{fig::proj_mzdmd}
\end{figure}

The difference in the predicted dynamics indicates that DMD uses only the resolved variables and cannot track the effect of the unresolved variables on resolved ones. 
Therefore, to take this effect into account, we implement the procedure proposed in~\cite{curtis2021dynamic}.
However, instead of constructing the first-order approximation of the optimality condition, we use the automatic differentiation technique~\cite{baydin2018automatic,jax2018github} to solve the stated optimization problem~\eqref{eq::mzdmd_final_opt} explicitly with the adaptive gradient method Adam~\cite{kingma2014adam}.
The parameters of the MZ-DMD approach are the following: learning rate is $10^{-3}$, number of iterations for internal optimization with Adam optimizer is $5$.
The initial value for the unknown matrix in MZ-DMD is equal to $\bA_{dmd}$.
Denote by $N_u$ the number of vectors $\bn \sim \mathcal{N}(0, \sigma^2\bI)$ used to initialize the memory term.
We use $N_u = 100$ random vectors for the initialization of vector  $\bn$.
To generate the result dynamics, we average the spectrum of the obtained matrices for every $\bn$ over $N_u$ generated samples and use equation~\eqref{eq::dynamic_recontr}, where $(\hat{x}_1, \hat{x}_2) = (1, 0)$. 
We also evaluate the empirical variance of the obtained expected dynamics following~\eqref{eq::variance_estimate}.
The results of the comparison of projected dynamics and dynamics induced by MZ-DMD  approach are presented in Figure~\ref{fig::proj_mzdmd}.


\begin{figure}[!ht]
    \centering
    \subfloat{\includegraphics[width=0.45\textwidth]{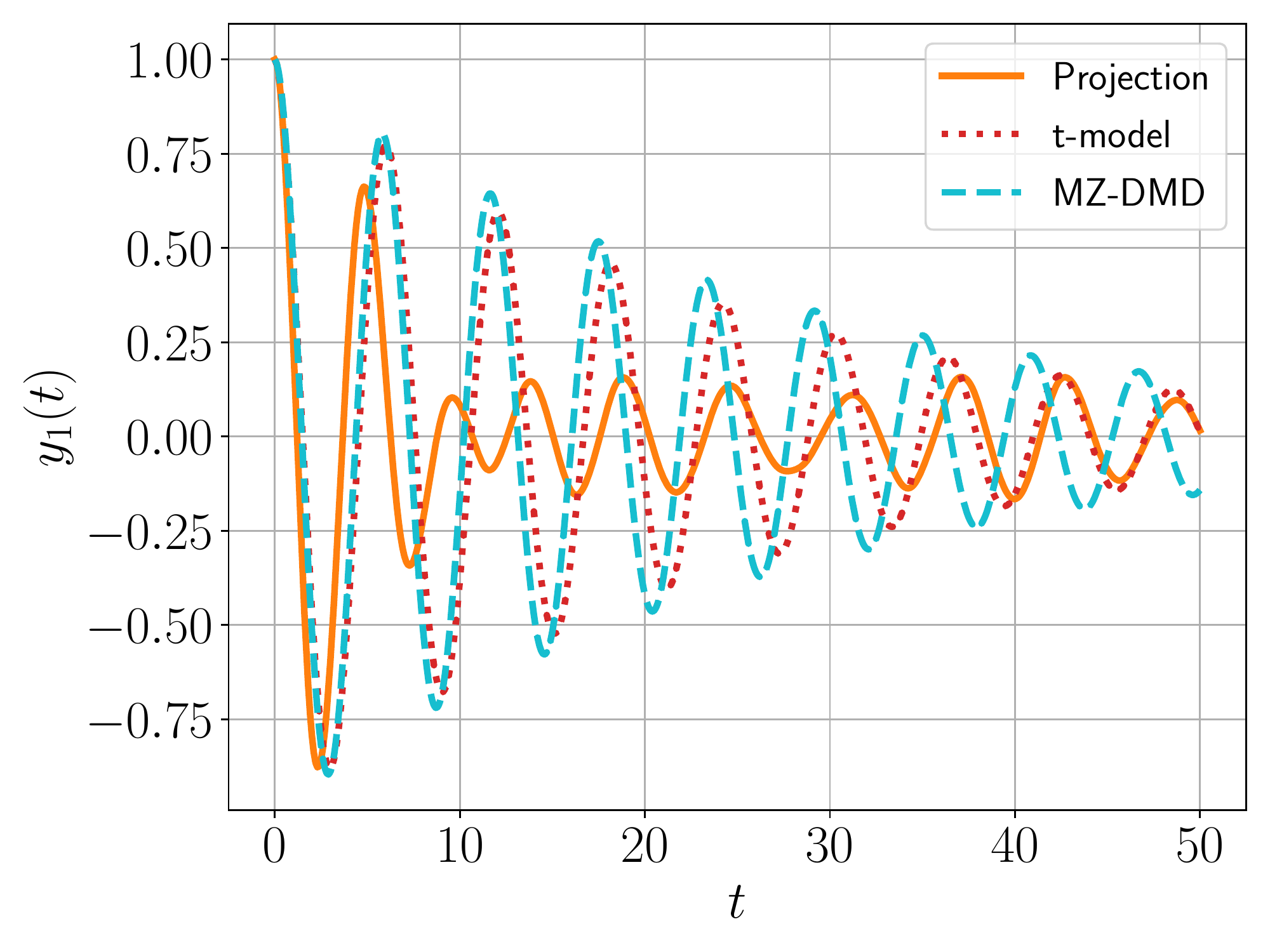}}
    ~
    \subfloat{\includegraphics[width=0.45\textwidth]{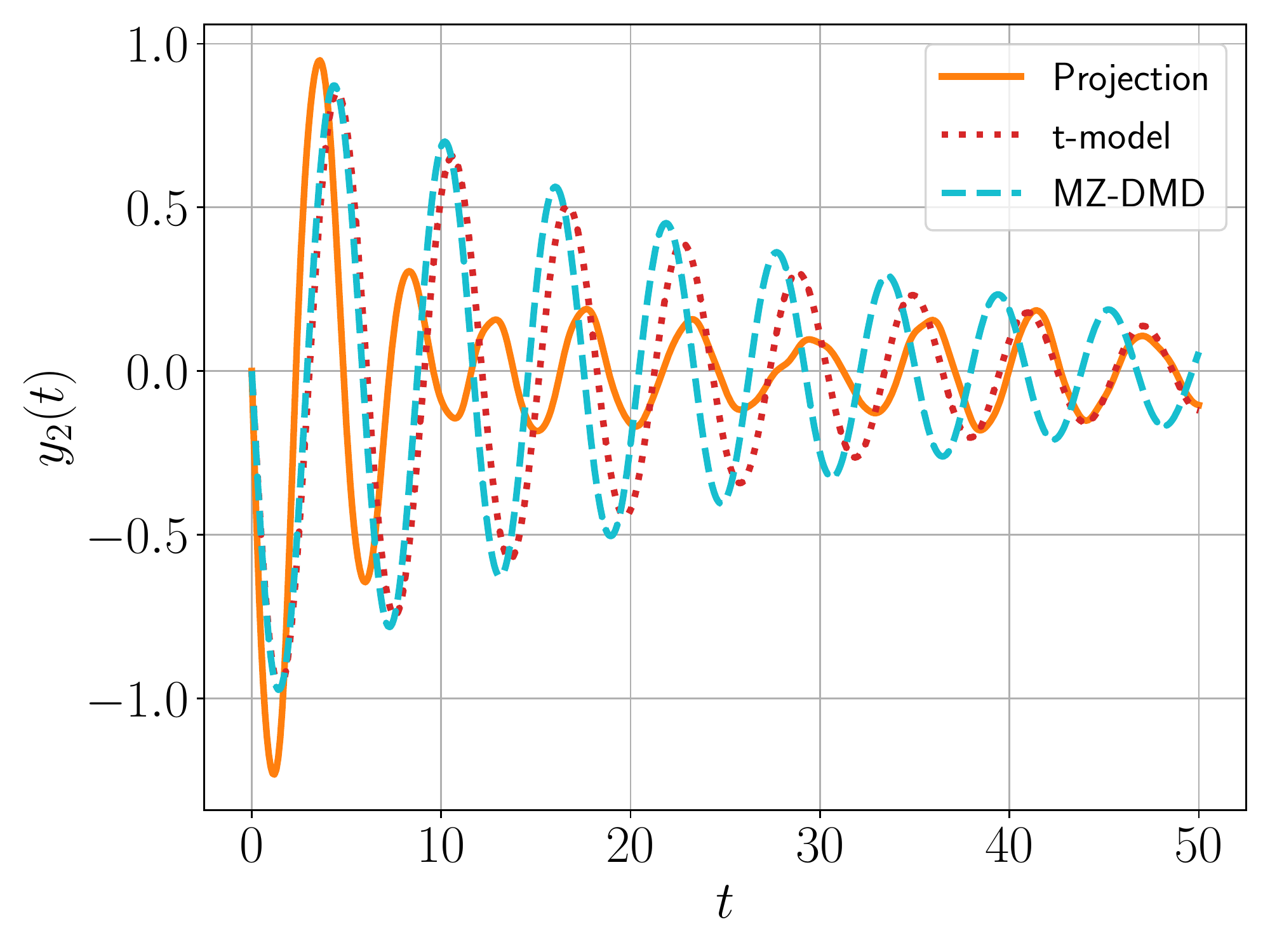}}
    \caption{Comparison of the projected dynamics, MZ-DMD reconstruction and the t-model for resolved $y_1(t)$ and $y_2(t)$. Here, we see that the t-model provides the similar dynamics as MZ-DMD approach and even better approximates Projection dynamics although requires less computationally intensive operations. We use $N_u = 100$ initializations of vectors $\bn$. The variance of t-model estimate is smaller than $10^{-1}$ uniformly on the considered time interval and is not shown in plot.}
    \label{fig::tmodel_mzdmd}
\end{figure}

The use of the t-model is illustrated in Figure~\ref{fig::tmodel_mzdmd}.
This modification accelerates solving the optimization problem since the gradient can be computed much faster.
Moreover, we observe in Figure~\ref{fig::tmodel_mzdmd} that the t-model provides a similar dynamics to the MZ-DMD approach and is sufficiently accurate to model the system described by~\eqref{eq::simple_model}.
We use the same averaging technique over $N_u$ randomly sampled vectors $\bn$ that is used to obtain the dynamics for MZ-DMD approach. 


\paragraph{Running time comparison} 
The aforementioned experiments were carried out on a single NVIDIA Tesla V100 GPU with JAX framework.
The running time of the solving optimization problems for multiple initializations of $\bn$ are 46 and 30 seconds for MZ-DMD and t-model, respectively.
The result of the Projection method is computed for 10 minutes.
To accelerate computations, we use the Just-In-Time compilation technique and apply it to gradient computing functions.
From the theoretical point of view, the complexity of gradient computation with the  automatic differentiation is bounded from above with some constant multiplied by the complexity of the function computation~\cite{griewank2008evaluating}.
Therefore, we have some theoretical guarantee that the optimization process requires a limited amount of time for more complex dynamical systems with unresolved variables.

\paragraph{Apriori knowledge of the model property} 
We observe that a limitation of the proposed method is the apriori knowledge about coupling  the variables.
For example, if the resolved and unresolved variables are decoupled, then the memory term has to be initialized with zero vector.
Otherwise, one can get incorrect expected dynamics that contradict the property of the original dynamical system.
This requirement makes MZ-DMD and t-model approaches model specific.
In future work, we are going to address this issue and investigate modifications of the proposed approach that will automatically adjust the initialization of the memory term.

\section{Conclusion}
\label{sec::conclusion}

In this study we develop a method to analyze incompletely measured data such that one variables are resolved and other ones are unresolved.
We demonstrate that the standard DMD approach and its non-linear modifications can not capture the expected dynamics since they assume complete description of samples.
To model the dynamics of the resolved variables, we consider the  Mori-Zwanzig representation and approximate it with so-called t-model.
From the t-model we state the optimization problem with respect to the transition operator based on the resolved variables.
Our method generalizes the DMD method for the incomplete measurement setting.
The proposed method requires essentially less computational resources than the previously proposed combination of exact optimal prediction and DMD method.
We present the simulation results for a test Hamiltonian system, which confirm that our method establishes sufficiently accurate approximation of the expected dynamics of resolved variables with limited variance.
The approach allows us to obtain a reasonably good prediction of the averaged dynamics from only one measurement with incomplete data. 
The predicted dynamics corresponds to the conditional expectation under available incomplete information.

\section*{Acknowledgment}
The work presented in Sections 1 and 2 was carried out as a part of the AMPaC Megagrant project supported by Skoltech and The Ministry of Education and Science of Russian Federation, Grant Agreement No 075-10-2021-067, Grant identification code 000000S707521QJX0002.
Section 3 was supported by the grant EPSRC EP/V038249/1.
Section 4 was supported by Ministry of Science and Higher Education grant No. 075-10-2021-068.

\bibliographystyle{unsrt}

\end{document}